\documentstyle{article}

\newcommand{\Frac}[2]{\displaystyle\frac{#1}{#2}}
\newcommand{\Sum}{\displaystyle\sum}
\newcommand{\Prod}{\displaystyle\prod }
\newcommand{\N}{\mbox{\rm I\kern-1.5pt N}}
\newcommand{\R}{\mbox{\rm I\kern-1.5pt R}}
\newcommand{\nsubseteq}{\subseteq\!\!\!\!\!\!/}
\newcommand{\D}{\mbox{${\rm I} \! {\rm D}$}}
\newcommand{\F}{\mbox{${\rm I} \! {\rm F}$}}
\newcommand{\PP}{\mbox{${\rm I} \! {\rm P}$}}
\newcommand{\Z}{\mbox{\sf Z \hspace{-1.1em} Z}}
\newcommand{\dem}{\par{\em Proof\/: }\\ \noindent }
\newcommand{\findemo}{$\;\;\Box$\\}
\newfont{\frack}{eufm10}

\newtheorem{defi}{Definition}[section]
\newtheorem{prop}[defi]{Proposition}
\newtheorem{thm}[defi]{Theorem}
\newtheorem{nota}[defi]{Remark}
\newtheorem{algor}[defi]{Algorithm}
\newtheorem{ejplo}[defi]{Example}

\begin{document}

\title{Symbolic Hamburger-Noether expressions of plane curves 
and construction of AG codes}
\author{A. Campillo$\;\,^{\ast}$ \and J. I. Farr\'{a}n 
\thanks{Both authors are partially supported by DIGICYT PB97-0471. 
A. Campillo is with Dpto. Algebra, Univ. Valladolid, Spain, e-mail: campillo@agt.uva.es, and 
J. I. Farr\'{a}n is with Dpto. Matem\'{a}tica Aplicada a la Ingenier\'{\i}a, Univ. Valladolid, Spain, e-mail: ignfar@eis.uva.es }}
\date{September 15, 1999}
\maketitle

\begin{abstract}

We present an algorithm to compute bases for the spaces ${\cal L}(G)$ 
and $\Omega(G)$\/, provided $G$ is a rational divisor over 
a non-singular absolutely irreducible algebraic curve, 
and also another algorithm to compute the Weierstrass semigroup at $P$ 
together with functions for each value in this semigroup, provided $P$ 
is a rational branch of a singular plane model for the curve. 
The method is founded on the Brill-Noether algorithm by combining 
in a suitable way the theory of Hamburger-Noether expansions and 
the imposition of virtual passing conditions. Such algorithms are 
given in terms of symbolic computation by introducing the notion of 
symbolic Hamburger-Noether expressions. 
Everything can be applied to the effective construction of Algebraic 
Geometry codes and also in the decoding problem of such codes, 
including the case of the Feng and Rao scheme for one-point codes. 

{\bf MSC2000}: 14 Q 05 -- 11 T 71


\end{abstract}

\section{Introduction}

Since the construction at the beginning of the 80's by Goppa 
of linear codes using Algebraic Geometry (see \cite{Goppa}), 
the theory of algebraic-geometric codes has been extensively 
developed. Algebraic Geometry codes (AG codes in short) can 
be constructed from any smooth algebraic projective curve 
$\tilde{\chi}$ defined over a finite field $\F$ as images of 
linear maps involving either residues at certain rational points 
of differential forms in $\Omega(G-D)$\/, or evaluations at such 
points of rational functions in ${\cal L}(G)$\/, where $G$ is 
a rational divisor over $\F$ and $D=P_{1}+\ldots+P_{n}$ is the 
formal sum of the rational points above considered, $G$ and $D$ 
having disjoint supports (details in section {\bf 6}). 

In order to construct such codes, the 
main difficulty in practice is the computation of vector bases 
for the spaces ${\cal L}(G)$\/, since the construction with 
differential forms is derived by duality. 
When the curve is given by means of a plane singular birational 
model $\chi$\/, some general methods can be used for this task 
if one knows well enough the singularities, namely the 
Brill-Noether \cite{HachTh} and Coates methods \cite{DuvTh}. 

On the other hand, nice codes need to have also good decoding 
algorithms. From the beginning of the 90's several decoding 
methods have been developed (see \cite{HohPel} for a survey on this matter). 
In the case $G=mP$ for some extra rational point $P$ and $m>0$, 
Feng and Rao gave in \cite{FR} a simple method based on a majority 
voting test, nowadays considered to be the most efficient decoding 
procedure. This method requires the previous knowledge of the 
Weierstrass semigroup of $\chi$ at the rational branch given by $P$\/, 
together with a rational function $f_{l}\in\F(\chi)$ regular outside 
$P$ and achieving a pole at $P$ of order $l$\/, for each $l$ in this 
semigroup, as we will show in section {\bf 6}. 

Again, the main difficulty turns out to be the computation of 
Weierstrass semigroups and such functions $f_{l}\,$. 
This difficulty is in fact the main obstacle for practical uses in 
Coding Theory, among others, of the construction by Garc\'{\i}a and 
Stichtenoth in \cite{GStich} of a sequence of curves achieving 
the Drinfeld-Vl\u{a}du\c{t} bound given in \cite{D/V}. 

Thus, the effective coding and decoding of AG codes depend on the 
resolution of two basic problems: computation of a vector basis 
for ${\cal L}(G)$ and the computation of Weierstrass semigroups 
together with functions achieving its values. 
The objective of this paper is to give a complete symbolic-computation 
treatment of these basic problems from the knowledge of a singular plane 
birational model $\chi$ for the smooth curve $\tilde{\chi}$\/, 
what is actually the most usual way to give a curve. 
The algorithms of this paper are at present being implemented 
by the authors in the computer algebra system SINGULAR \cite{Sing}, 
created by Greuel, Pfister and Schoenemann. 

Our approach is based on very classical ideas. First we consider 
Hamburger-Noether expansions from the symbolic viewpoint; 
more precisely, we introduce in the paper the so-called 
symbolic Hamburger-Noether expressions, which will provide us 
with both all the information on the singularities and (symbolic) 
parametrizations for all their rational branches. 
Hamburger-Noether expansions are developed in \cite{Camp} 
for the case of irreducible curve singularities over algebraically 
closed fields. Here we will need not only the symbolic version but 
also the case of general plane curve singularities over perfect fields 
(finite fields in practice) as developed in \cite{ACJC}. 

In particular, from the knowledge of the singularities one can compute 
the adjunction divisor, and from it the imposing conditions test for 
being an adjoint (see sections {\bf 3} and {\bf 4}). 
This becomes important for the approach to the first basic problem, 
since via the Brill-Noether algorithm 
the computation of a vector basis for ${\cal L}(G)$ is reduced 
to computing vector bases for some concrete spaces of adjoints, 
which are obtained by imposing certain assigned conditions. 
We show how the computation of such adjoint bases can be done by using 
the so-called {\em principle of discharge} due to Enriques in \cite{Enri} 
(see \cite{Casas} and \cite{Lip} for a modern treatment). 
Thus, our solution to the first basic problem is derived from the three 
classical theories of Hamburger-Noether (section {\bf 2}), 
Brill-Noether (section {\bf 3}) and Enriques (section {\bf 4}). 

The second basic problem is approached in similar terms. 
In fact, from a singular plane model the adjunction theory of plane curves 
can be applied to give an algorithm to compute the Weierstrass semigroup 
and the corresponding functions (see section {\bf 5}). 
Again, this algorithm becomes effective using symbolic 
Hamburger-Noether expressions at the singularities of $\chi$\/. 
Finally, we apply these methods to the construction of AG codes 
in section {\bf 6}.

\section{Symbolic Hamburger-Noether expressions 
of plane curve singularities}

In this section, we will introduce the symbolic Hamburger-Noether expressions 
for a plane curve singularity. For this, we fix in the sequel 
an arbitrary perfect field $\F$ and an absolutely 
irreducible projective algebraic plane curve $\chi$ defined over $\F$. 
For a closed point $P$ of $\chi$ 
with local ring $R={\cal O}_{\chi,P}$ we denote by a rational branch of $\chi$ 
at $P$ any maximal prime ideal of $\overline{R}$\/, where $\overline{R}$ 
denotes the the semilocal ring given by the normalization of $R$\/. 
The datum of such a maximal ideal is equivalent to give a minimal prime ideal 
of $\widehat{R}$\/, the completion with respect to the 
Jacobson radical of $R$ (see \cite{Camp} and \cite{ACJC} for details).

Assume that we have chosen an affine chart containing $P$\/, and let 
$A=\F[{\rm X},{\rm Y}]/(f({\rm X},{\rm Y}))$ be the affine ring of 
coordinates, $f({\rm X},{\rm Y})=0$ being the affine equation of 
the curve in this chart. Regarding $P$ as a non-zero prime ideal 
of $A$\/, one has $k(P)\cong K\hookrightarrow\widehat{A_{P}}$ and, 
for practical reasons, one can actually write 
$K=\F[{\rm Z}]/(Q({\rm Z}))$ for an irreducible polynomial $Q\in\F[{\rm Z}]$. 
Thus, up to a translation in $K[{\rm X},{\rm Y}]\,$, we can assume that $P$ 
is the origin, the defining ideal of $P$ being then $({\rm X},{\rm Y})$\/. 

With these notations, one has 
$\widehat{R}\cong K[[{\rm X},{\rm Y}]]/(f({\rm X},{\rm Y}))$\/, 
and hence there exists a natural morphism 
$K[[{\rm X},{\rm Y}]]\rightarrow\widehat{R}$. 
This allows us to introduce the following definition.

\begin{defi}

In the above conditions, a rational parametrization of $\chi$ at $P$ 
related to the coordinates ${\rm X},{\rm Y}$ is a $K$\/-algebra morphism 
$$r\;:\;K[[{\rm X},{\rm Y}]]\rightarrow K_{1}[[t]]$$
continuous for the $({\rm X},{\rm Y})$\/-adic and the $t$\/-adic 
topologies, such that $Im(r)\nsubseteq\,K_{1}$ and $f\in ker(r)$\/, 
where $K_{1}$ is a finite extension of $K$ and $t$ is an indeterminate. 
This is equivalent to give formal series $x(t),y(t)\in K_{1}[[t]]$ 
with at least one non identically zero such that $f(x(t),y(t))\equiv 0$. 

\end{defi}

We can associate to each rational parametrization $r$ the rational branch 
given by the minimal prime ideal $\mbox{\frack p}=ker(\widehat{r})$, where 
$\widehat{r}\;:\;\widehat{R}\rightarrow K_{1}[[t]]$ is the natural morphism 
induced by $r$\/. Thus, we say that $r$ is a rational parametrization of 
the branch $\mbox{\frack p}$\/. 

We say that another rational parametrization 
$s:K[[{\rm X},{\rm Y}]]\rightarrow K_{2}[[u]]$ is derived from $r$, 
and it is denoted by $s\succ r$, if there exists a formal series 
$t(u)\in K_{2}[[u]]$ with positive order and a $K$\/-algebra morphism 
$\sigma\;:\;K_{1}[[t]]\rightarrow K_{2}[[u]]$ with 
$\sigma(t)=t(u)$, such that $s=\sigma\circ r$\/. 
One has that $\succ$ is a partial preorder, and we say that two rational 
parametrizations $r$ and $s$ are equivalent if $s\succ r$ and $r\succ s$. 
Thus, a rational parametrization $r$ is called {\em primitive} if it is 
minimal (with respect to the partial preorder $\succ$\/) modulo equivalence, 
and moreover the extension $K_{1}|\F$ is also minimal (that is, 
$r({\rm X})$ and $r({\rm Y})$ are not both in $K'[[t]]$ for some field $K'$ 
with $\F\subseteq K'\subset K_{1}$ and $K'\neq K_{1}$\/). One actually has that 
rational branches at $P$ are in bijection with equivalence classes of 
primitive rational parametrizations at $P$ and, in particular, there always exist 
rational parametrizations (details again in \cite{ACJC}). By choosing a 
primitive rational parametrization for each rational branch one obtains 
a so-called {\em standard set of rational parametrizations} at $P$\/, 
and our next aim is the effective computation of such a set 
by means of the so-called Hamburger-Noether expansions.

Although a general definition can be given for arbitrary singular 
curves, we will study only the case of plane curves in order to get 
effective computations. 
Thus, let $\rho\;:\;\F[[{\rm X},{\rm Y}]]\rightarrow F[[u]]$ be 
a rational parametrization over $\F$ of the plane curve $\chi$ defined 
at the point $P$\/, $F$ being a finite extension of the base field $\F$. 
Denote in short by ${\cal O}$ the local ring of $\chi$ at $P$\/. 
One can consider ${\cal O}$ in fact as a subring of 
$\F[[u]]$, the images of ${\rm X}$ and ${\rm Y}$ being 
a minimal system of generators of the maximal ideal of ${\cal O}$\/. 

\begin{defi}

We introduce the $\,${\bf Hamburger-Noether expansion} of $\chi$ at $P$ 
for the branch given by $\rho$ to be a finite sequence $\D$ of 
expressions in the variables $Z_{-1},Z_{0},Z_{1},\ldots,Z_{r}$ of the form 
$$Z_{-1}=a_{0,1}Z_{0}+a_{0,2}Z_{0}^{2}+\ldots+a_{0,h_{0}}Z_{0}^{h_{0}}+
Z_{0}^{h_{0}}Z_{1}$$
$$Z_{0}=a_{1,2}Z_{1}^{2}+a_{1,3}Z_{1}^{3}+\ldots+a_{1,h_{1}}Z_{1}^{h_{1}}+
Z_{1}^{h_{1}}Z_{2}$$
$$..........................................................$$
$$Z_{r-2}=a_{r-1,2}Z_{r-1}^{2}+a_{r-1,3}Z_{r-1}^{3}+\ldots+
a_{r-1,h_{r-1}}Z_{r-1}^{h_{r-1}}+Z_{r-1}^{h_{r-1}}Z_{r}$$
$$Z_{r-1}=\sum_{i\geq1}a_{r,i}Z_{r}^{i}$$
where $r$ is a non-negative integer, $a_{j,i}\in F$, $a_{k,1}=0$ if $k>0$, 
$h_{j}$ are positive integers and moreover 
$$f(Z_{0}(Z_{r}),Z_{-1}(Z_{r}))=0\;\;\;{\rm in}\;\;\;F[[Z_{r}]]$$
$f\in\F[[{\rm X},{\rm Y}]]$ being a generator of the ideal $ker\,(\rho)\,$\/. 

\end{defi}

The existence of such expansions and the finiteness of the number of 
lines is refered to \cite{Camp}, \cite{ACJC} or \cite{HN}. 
In fact such an expansion $\D$ 
always gives a primitive rational parametrization equivalent to 
$\rho$ if we consider ${\rm X}\equiv{\rm Z}_{0}$ and 
${\rm Y}\equiv{\rm Z}_{-1}$ as a function of the local parameter 
$s=Z_{r}$ by successive substitutions. Moreover, $\D$ only depends 
on the branch given by $\rho$ and the choice of the parameters 
$x,y$ in ${\cal O}$ given by the images of ${\rm X},{\rm Y}$ under $\rho$\/. 
Thus, for ${\rm X}$ and ${\rm Y}$ fixed the (finite) set of all the 
possible non-equivalent Hamburger-Noether expansions form a standard 
set of rational parametrizations of $\chi$ at $P$ (see \cite{ACJC}).  

\begin{nota}

The role played by the Hamburger-Noether expansions in arbitrary 
characteristic is just the same as that classically played by the 
Puiseux expansions in characteristic $0$, which are given by 
$$X(t)=\alpha t^{\nu}$$
$$Y(t)=\sum_{i\geq\nu}\lambda_{i}t^{i}$$
where $\alpha\in F^{\ast}$ and $\lambda_{i}\in F$\/. 
The main problem of the Puiseux expansions is that they do not always exist 
in positive characteristic, and when such expansions exist they are 
rational but the problem of making them primitive is not at all trivial 
(see \cite{Camp} or \cite{DuvTh}). These are the reasons why we use 
Hamburger-Noether expansions. 

\end{nota}

Now we show how to compute the Hamburger-Noether expansions without 
having a priori any local parametrization of the branch, but only with 
the aid of the Newton polygon of the local equation of $\chi$ at $P$\/. 
We will do it for the case of only one rational branch at $P$ for the 
sake of simplicity, but the method also works for several branches 
(in the reduced case) because of the fact that the Newton polygon would 
be the collection of those of each branch joined together with increasing 
slope (see \cite{HN} for further details). 

More precisely, let $\F$ be a perfect field and let $\chi$ be given 
in affine coordinates by the local equation $f({\rm X},{\rm Y})=
\sum_{\alpha,\beta\geq 0}c_{\alpha\beta}{\rm X}^{\alpha}{\rm Y}^{\beta}=0$, 
$f$ being an irreducible polynomial in $\F[{\rm X},{\rm Y}]$\/. 
Assume that we want to study the point $P=(0,0)$ and that there is 
only one rational branch at the origin defined over $\F$. 
Then we consider the Newton diagram of $f$ 
$$D(f)\doteq\{(\alpha,\beta)\;|\;c_{\alpha\beta}\neq 0\}$$
and we call {\em Newton polygon} of $f$ (at the origin) the set of 
all the bounded segments of the convex hull of 
$D(f)+\R^{2}_{+}\;$, and it will be denoted by $P(f)$\/. 

Excluding the trivial cases where the curve is one of the coordinate 
axes, let $l$ (respectively $n$\/) be the minimum integer such that 
$(l,0)\in D(f)$ (respectively $(0,n)\in D(f)$\/). We can obviously 
assume that $n\leq l$. In this case, the Newton polygon consists just of 
one segment with non-zero slope and extremes $(l,0)$ and $(0,n)$\/. 

If $\Delta=P(f)$ is the Newton polygon we can define 
$$L({\rm X},{\rm Y})\doteq
\sum_{(\alpha,\beta)\in\Delta}c_{\alpha\beta}{\rm X}^{\alpha}{\rm Y}^{\beta}$$
One obviously has 
$L({\rm X},{\rm Y})=c\,D({\rm X},{\rm Y})$ for some $c\in\F^{\ast}$ 
and some $D({\rm X},{\rm Y})$ which is monic in ${\rm Y}$ 
and defined over $\F$\/. Moreover, by using the Hensel lemma one has 
$$D({\rm X},{\rm Y})=\Prod_{j=1}^{d}({\rm Y}^{n'}-\delta_{j}{\rm X}^{l'})^{e}$$
for some $\delta_{j}\in\overline{\F}^{\ast}$\/, where $ed=g\,c\,d\,(l,n)$\/. 
Then the {\em characteristic polynomial} of $\Delta$ is given by 
$$\Phi_{\Delta}(\lambda)\doteq\Prod_{j=1}^{d}(\lambda-\delta_{j})$$
It is an irreducible polynomial over $\F$ (that is, $\delta_{j}$ are 
conjugate each other by the Galois group over $\F$\/). Moreover, 
one has $l=l'ed$ and $n=n'ed$\/, being $g\,c\,d\,(l',n')=1$. 

If we write $l=qn+h$ with $0\leq h<n$, we find one of the following two cases: 

\begin{description}

\item[Case 1: ] $h=0$\/, what implies $ed=n$, $l'=q$ and $n'=1$. 
Thus write  
$$a_{0,1}=\ldots=a_{0,l'-1}=0,\;\;\;{\rm and}\;\;\;a_{0,l'}=\delta$$
$\delta$ being a symbolic root 
\footnote{We mean by a symbolic root of $\Phi_{\Delta}(\lambda)$ that one 
substitutes $\F$ by the field $\F_{1}=\F[\lambda]/(\Phi_{\Delta}(\lambda))$ 
and one takes as $\delta$ the residual class of $\lambda$ in this field. } 
of $\Phi_{\Delta}(\lambda)$\/, 
we get that the first line of the Hamburger-Noether expansion starts with 
$$Z_{-1}=a_{0,l'}Z_{0}^{l'}+\ldots$$
Then we transform $f$ by 
$$T_{1}(f,\delta,l')=f({\rm X},{\rm Y}+\delta\,{\rm X}^{l'})=
f_{1}({\rm X},{\rm Y})$$
getting $f_{1}$ with a segment of extremes $(l_{1},0)$ and $(0,n)$ 
as the Newton polygon, being $l_{1}>l$, and we iterate the process\footnote{Notice 
that this process could terminate if there is no point of the form $(l_{1},0)$\/. }, 
taking into account that $f_{1}$ has the coefficients in the field 
$\F_{1}=\F[\lambda]/(\Phi_{\Delta}(\lambda))$ 
and that it is irreducible over such field. 

\item[Case 2: ] $h>0$\/; in this case, the first line of the 
Hamburger-Noether expansion is just 
$$Z_{-1}=Z_{0}^{q}Z_{1}$$
Now, since the polynomial 
$U(f,l,n)=f({\rm Y},{\rm X}\,{\rm Y}^{q})$ is divisible by ${\rm Y}^{nq}$, 
we can transform $f$ by 
$$T(f,l,n)=\Frac{f({\rm Y},{\rm X}\,{\rm Y}^{q})}{{\rm Y}^{nq}}=
f_{1}({\rm X},{\rm Y})$$
Thus the obtained Newton polygon $\Delta_{1}$ has 
$(n,0)$ and $(0,h)$ as extremes, being $h<n$, and its characteristic 
polynomial is 
$\Phi_{\Delta_{1}}(\lambda)=\lambda^{e}\,\Phi_{\Delta}(1/\lambda)$\/. 
Then we repeat the process, looking for the next line of the 
Hamburger-Noether expansion, identifying in $T(f,l,n)$ 
${\rm X}\equiv{\rm Z}_{1}\,$ and $\,{\rm Y}\equiv{\rm Z}_{0}\,$. 

\end{description}

If the case $h=0$ is found some consecutive times during the computation 
of the line $k+1$ of the Hamburger-Noether expansion, where 
$k\in\{0,1,\ldots,r-1\}$\/, 
we append the new obtained result to the previous part of that line  
until we get the case $h>0$, where we append the term $Z_{k}^{q}Z_{k+1}$ 
and change to the next line. 

Hence, by applying a finite number of transformations of type $T_{1}$ or $T$, 
we get a trivial Newton polygon where either $\min\,(l,n)=1$ or there is 
no point in the vertical axis, and the procedure stops. In this case 
one has all the lines of the Hamburger-Noether expansion except the 
last one. But then the equation $f({\rm X},{\rm Y})$ is transformed 
into a polynomial $g({\rm Z}_{r},{\rm Z}_{r-1})$\/, $g$ being defined 
over a field $\F'$ which is obtained by successive symbolic extensions 
of $\F$ and one has $\Frac{\partial g}{\partial{\rm Z}_{r-1}}(0,0)\neq 0$. 
This means that one can obtain as many terms as needed of the last line 
of the Hamburger-Noether expansion from the polynomial 
$g({\rm Z}_{r},{\rm Z}_{r-1})$ as an implicit function 
(that is, by indeterminate coefficients), since this line 
represents ${\rm Z}_{r-1}$ as a formal series in the variable ${\rm Z}_{r}\,$. 
Thus, we do not need in practice the data given by the (infinite) series 
of the last line of the Hamburger-Noether expansion, but only the (finite) 
data of the implicit equation $g$\/, which contains the same information. 
This data is what we call a {\em symbolic Hamburger-Noether expression}\/, 
and it can be computed in an effective way by the above method for every 
singular closed point of $\chi$ (initially written in a symbolic extension 
of the base field if such a point is not rational). Even more, we do 
all these computations in successive symbolic extensions of $\F$ instead 
of considering a sufficiently large extension of it, 
what in practice saves a lot of time.

In the case of several branches the characteristic polynomial is not 
irreducible and each branch corresponds to an irreducible factor 
of this polynomial and its corresponding symbolic root, proceeding 
as in the case of one branch with every factor in parallel. Hence, 
in the general case we have to add in each step of the previous 
algorithm a factorization procedure for the corresponding characteristic 
polynomial, what also has an effective solution. 
Each irreducible factor follows at least one of the rational branches, 
so that one has an algorithm in form of tree. Thus, the branches of the 
tree given by this algorithm correspond one-to-one to the branches of 
the curve at the considered point, and for each tree branch one has 
associated (as a byproduct of the algorithm) the symbolic 
Hamburger-Noether expression corresponding to the curve branch. 
The computation of Hamburger-Noether expansions is a known method 
and it has been implemented with the computer algebra system 
SINGULAR \cite{Sing}.

\begin{ejplo}

Let $\chi$ be the projective plane curve over $\F_{2}$ given by 
$$F({\rm X},{\rm Y},{\rm Z})={\rm X}^{10}+{\rm Y}^{8}{\rm Z}^{2}+
{\rm X}^{3}{\rm Z}^{7}+{\rm Y}{\rm Z}^{9}=0$$
with the only singular point $P=(0:1:0)$ which is rational over $\F_{2}\,$, 
being furthermore the unique point of $\chi$ at infinity. 
Take the local equation 
$$f({\rm X},{\rm Z})={\rm X}^{10}+
{\rm X}^{3}{\rm Z}^{7}+{\rm Z}^{9}+{\rm Z}^{2}$$ 
of $\chi$ where $P$ is the origin, and apply the Hamburger-Noether algorithm 
to this equation. 

With the above notations, one has 
$L({\rm X},{\rm Z})={\rm Z}^{2}+{\rm X}^{10}=({\rm Z}+{\rm X}^{5})^{2}\,$; 
thus $l=10$, $l'=5$, $n=e=2$, $n'=d=1$ and $q=5$, being in the case $h=0$. 

The characteristic polynomial is $\Phi(\lambda)=\lambda+1$ 
and thus the symbolic root is nothing but $\delta=1$, that is, 
we do not need to enlarge the base field $\F_{2}\,$. Hence, one has 
$$a_{0,0}=\ldots=a_{0,4}=0\;\;\;\;\;\;\;\;a_{0,5}=1$$
and we do the change 
$$f_{1}({\rm X},{\rm Z})=f({\rm X},{\rm Z}+{\rm X}^{5})=
{\rm Z}^{2}+{\rm X}^{38}+\ldots$$
being now $L({\rm X},{\rm Z})=({\rm Z}+{\rm X}^{19})^{2}$ and thus 
$l=38$, $l'=19$, $n=e=2$, $n'=d=1$, $q=19$ and again $h=0$; 
one also has $\Phi(\lambda)=\lambda+1$ and $\delta=1$. Thus 
$$a_{0,6}=\ldots=a_{0,18}=0\;\;\;\;\;\;\;\;a_{0,19}=1$$
and we do the transform 
$$f_{2}({\rm X},{\rm Z})=f_{1}({\rm X},{\rm Z}+{\rm X}^{19})=
{\rm Z}^{2}+{\rm X}^{45}+\ldots$$
In this case, one has $L({\rm X},{\rm Z})={\rm Z}^{2}+{\rm X}^{45}$\/, 
obtaining $l=l'=45$, $n=n'=2$, $d=e=1$ and $q=22$, being now in the case 
$h=1>0$ and we have to change the line in the Hamburger-Noether expansion 
without enlarging the base field. 
Now the transform to do is 
$$f_{3}({\rm X},{\rm Z})=\Frac{f_{2}({\rm Z},{\rm X}{\rm Z}^{22})}{{\rm Z}^{44}}=
{\rm Z}+{\rm X}^{2}+\ldots$$
being now the origin a non-singular point of the new equation and the procedure 
ends with $r=1$. Thus, the symbolic Hamburger-Noether expressions at $P$ are 
$$\left\{\begin{array}{rcl}
Z_{-1}&=&Z_{0}^{5}+Z_{0}^{19}+Z_{0}^{22}Z_{1}\\
g(Z_{1},Z_{0})&=&
Z_{1}^{9}Z_{0}^{154}+Z_{1}^{8}Z_{0}^{151}+
Z_{1}^{8}Z_{0}^{137}+Z_{1}Z_{0}^{130}+
Z_{0}^{127}+Z_{1}^{7}Z_{0}^{113}+\\
\,&+&
Z_{1}^{6}Z_{0}^{110}+Z_{0}^{113}+
Z_{1}^{5}Z_{0}^{107}+Z_{1}^{4}Z_{0}^{104}+
Z_{1}^{3}Z_{0}^{101}+Z_{1}^{6}Z_{0}^{96}+\\
\,&+&
Z_{1}^{2}Z_{0}^{98}+Z_{1}Z_{0}^{95}+
Z_{1}^{4}Z_{0}^{90}+Z_{0}^{92}+
Z_{1}^{2}Z_{0}^{84}+Z_{1}^{5}Z_{0}^{79}+\\
\,&+&
Z_{1}^{4}Z_{0}^{76}+Z_{0}^{78}+Z_{1}Z_{0}^{67}+
Z_{1}^{4}Z_{0}^{62}+Z_{0}^{64}+Z_{0}^{50}+\\
\,&+&
Z_{1}^{3}Z_{0}^{45}+Z_{1}^{2}Z_{0}^{42}+
Z_{1}Z_{0}^{39}+Z_{0}^{36}+
Z_{1}^{2}Z_{0}^{28}+Z_{0}^{22}+\\
\,&+&
Z_{1}Z_{0}^{18}+Z_{0}^{15}+Z_{1}Z_{0}^{11}+
Z_{0}^{8}+Z_{1}^{2}+Z_{0}
\end{array}\right.$$

\vspace{.2cm}

\end{ejplo}

\section{Normalization, resolution and adjunction via 
symbolic Hamburger-Noether expressions}

The purpose of this section is the revision of some classical concepts 
taking into account the symbolic Hamburger-Noether expressions 
which have been introduced in the previous section. 
Thus, for a given plane curve $\chi$ one can consider its normalization, 
that is the proper birational morphism 
$${\bf n}\;:\;\tilde{\chi}\rightarrow\chi$$
where $\tilde{\chi}$ is the curve obtained by gluing together 
the affine charts given by the normalization of the affine 
graded $\F$\/-algebras $A_{U}$ for all affine charts $U$ of $\chi$\/. 
The curve $\tilde{\chi}$ can be obtained as the blowing-up of the 
conductor, that is the sheaf of ideals locally given by 
$${\cal C}_{\chi}(U)\doteq\{f\in\overline{{\cal O}_{\chi}(U)}\;|\;
f\,\overline{{\cal O}_{\chi}(U)}\subseteq{\cal O}_{\chi}(U)\}$$
Nevertheless, it is better in practice to look at $\tilde{\chi}$ as 
successive blowing-ups of all the closed points of $\chi$ which are 
singular until we get a curve without singular points, since this approach 
can be explicitly described by equations. 
In each of those blowing-ups one has as result the corresponding {\em strict transform} 
$\chi_{i}$ for $i\geq 0$ (starting from $\chi_{0}=\chi$\/), defined as usual 
(see for example \cite{Thes} or \cite{HachTh}). This process can be represented 
by a combinatorial object called the {\em resolution forest} 
${\cal T}_{\chi}\,$, consisting of one {\em weighted oriented tree} 
for each singular closed point of $\chi$\/, and which is constructed as follows:

\begin{description}

\item[1)] The vertices represent the successive points which are obtained 
by blowing up singular points of the successive strict transforms $\chi_{i}$ of $\chi$ 
until one gets a non-singular point at the end of each branch of the process. 
Two such vertices $p$ and $q$ of one tree corresponding to the points 
$P$ and $Q$ are connected by an edge from $p$ to $q$ if $Q$ is one of the 
points obtained by blowing-up $P$\/. 

\item[2)] On each edge $\stackrel{\longrightarrow}{pq}$ of the forest we put a weight 
$\rho_{pq}\doteq \left[k(Q):k(P)\right]$, where $k(P)$ and $k(Q)$ are 
the corresponding residual fields of the local rings 
${\cal O}_{\chi_{i},P}$ and ${\cal O}_{\chi_{i+1},Q}$\/. 

\item[3)] If $p$ is the root of the tree corresponding to the singular 
point $P$ of $\chi$\/, then we put on $p$ an initial weight $[k(P):\F]$\/. 
On all the other vertices of the forest we can assign two 
alternative weights which are equivalent if we know the weights 
on the edges. In both cases one assigns to $p$ a weight for each 
branch of the tree passing through $p$\/, where by a branch we 
denote any upper extremal point of the forest, and we say that 
such a branch $q$ passes through $p$ there is an oriented path 
from $p$ to $q$ in ${\cal T}_{\chi}$ 
\footnote{Notice that such branches are in bijection with the rational 
branches at $P$ of the corresponding curve obtained by blowing-up $P$\/, 
and also with the closed points over $P$ of the normalization. }. 
The two alternative weights on $p$ for each $q$ are the following: 

\begin{description}

\item[(I)] The {\em multiplicity} at $P$ of the rational branch 
$\mbox{\frack q}$ corresponding to $q$ computed in the corresponding 
curve $\chi_{P}$ obtained by blowing-up $\chi$\/, that is the 
multiplicity $e_{p,q}$ of the noetherian ring 
$\widehat{{\cal O}_{\chi_{P},P}}/\mbox{\frack q}$ of dimension $1$ 
(denoting here $\mbox{\frack q}$ the corresponding minimal prime ideal 
of $\widehat{{\cal O}_{\chi_{P},P}}$). 

\item[(II)] The {\em order} at $P$ of the rational branch $\mbox{\frack q}$, 
that is the number $m_{p,q}\doteq\min\,\{\upsilon_{Q}(f)\;|\;
f\in\mbox{\frack m}_{\chi_{P},P}\}$, where $\mbox{\frack m}_{\chi_{P},P}$ is 
the maximal ideal of the local ring ${\cal O}_{\chi_{P},P}$ and 
$\upsilon_{Q}$ denotes the normalized valuation (that is, 
with $\Z$ as group of values) corresponding to $Q$ 
regarded as a point of $\tilde{\chi}$\/. 
The equivalence between both weights is given by the formula 
$$m_{p,q}\left[k(Q):k(P)\right]=e_{p,q}$$

\end{description}

\end{description}

Notice that the order is actually the multiplicity of each of the conjugate  
geometric branches lying over $P$\/, considering $\chi$ to be defined 
over the algebraic closure $\overline{\F}$ of $\F$. By substituting 
$\F$ by $\overline{\F}$ one obtains another combinatorial object which 
is much more complex than the one above described and that has all weights 
on the edges equal to $1$ and hence $m_{p,q}=e_{p,q}\,$. This object can be 
reconstructed from the rational object ${\cal T}_{\chi}$ 
(this is shown in \cite{Thes}) 
and it does not show properly the structure of $\chi$ over $\F$, being thus 
${\cal T}_{\chi}$ a more precise invariant of the normalization.

We will show now that from the computation of symbolic 
Hamburger-Noether expressions one gets, as a byproduct, 
the desingularization of the curve 
(see \cite{ACJC} and \cite{Thes} for more details). 
In fact, for simplicity consider again the case of only one rational branch. 
Let $f\in\F[{\rm X},{\rm Y}]$ be a local equation of $\chi$ at $P$\/, 
supposed rational and $P=(0,0)$ in the affine coordinates ${\rm X},{\rm Y}$ 
(otherwise we consider an initial symbolic extension $\F'$ instead of $\F$). 
If we write $l=qn+h$ as in the previous section, 
then the first $q$ {\em infinitely near} points $P=P_{0},P_{1},\ldots,P_{q-1}$ 
are rational over $\F$, being $P_{i}=(0,0)$\/, for $0\leq i\leq q-1$, in the 
local affine coordinates $\{ {\rm X},\Frac{\rm Y}{{\rm X}^{i}}\}$ at $P_{i}$\/.

If $h=0$, then $P_{q}$ has the symbolic field 
$\F_{1}=\F[\lambda]/(\Phi_{\Delta}(\lambda))$ as residual field, 
being $P_{q}=(0,0)$ in the local affine coordinates related to $\F_{1}$ 
given by $\{ {\rm X},\Frac{\rm Y}{{\rm X}^{q}}-\delta\}$\/, $\delta$ being 
a symbolic root of the characteristic polynomial $\Phi_{\Delta}(\lambda)$\/. 

If $h>0$, then the new coordinates are $\{ {\rm Z}_{1},{\rm Z}_{0} \}$\/, 
$P_{q}$ is rational over $\F$ and $P_{q}=(0,0)$ in these coordinates, 
${\rm Z}_{1}=0$ being now the exceptional divisor instead of ${\rm Z}_{0}=0$. 
Anyway, by doing successively the above changes of variables one easily gets 
the corresponding total, strict or virtual transform of any divisor. 

With this notation, the edges $\stackrel{\longrightarrow}{p_{i-1}p_{i}}$ of the resolution 
forest ${\cal T}_{\chi}$\/, $p_{j}$ corresponding to $P_{j}$\/, 
have weight $1$ either if $i<q$ or if $i=q$ and $h>0$, 
and weight $d$ if $i=q$ and $h=0$. The value $e\cdot n'$ in each step is just 
the order of that branch at $P_{0},\ldots,P_{q-1}\,$, 
and $n=d\cdot e\cdot n'$ is the multiplicity.  
The weights at $P_{q}$ appear in the next step of the algorithm, 
where $P_{q}$ plays the role of $P_{0}=P$\/, and so on.

When one gets the trivial polygon by iterating this method, one obtains all 
the infinitely near points with all the weights of the combinatorial object 
${\cal T}_{\chi}$\/. When the procedure ends, one has the coordinates 
$\{ {\rm Z}_{r},{\rm Z}_{r-1} \}$ and the local equation 
$g({\rm Z}_{r},{\rm Z}_{r-1})$\/, satisfying 
$\Frac{\partial g}{\partial{\rm Z}_{r-1}}(0,0)\neq 0$. Doing $s$ 
additional transformations of type $T_{1}$ one obtains the embedded resolution, 
being ${\rm Z}_{r}^{s}$ the initial form of $g({\rm Z}_{r},{\rm Z}_{r-1})$\/.

In the case of several branches, the resolution can be obtained taking into 
account that there are as many irreducible factors of the characteristic 
polynomial as infinitely near points in the exceptional divisor, and the 
corresponding symbolic roots yield suitable local coordinates for such 
points, that is, everything can be done, branch by branch, 
with an algorithm in form of a tree.

\begin{ejplo}

In the example 2.4 , 
one obtains the resolution tree of $\chi$ at $P$ as the sequence of points 
$$P\equiv p_{0}\rightarrow p_{1}\rightarrow\ldots\rightarrow p_{21}\rightarrow p_{22}\equiv q$$
corresponding to rational points of multiplicity $e_{p_{i},q}=2$ if 
$i=0,\ldots,21\,$, and $e_{p_{22},q}=1$, the weights of all the edges being $1$ 
as the initial weight, since we have never enlarged the base field. 

\end{ejplo}

A useful information which one can derive from ${\cal T}_{\chi}$ is 
the adjunction divisor ${\cal A}$ of the singular plane curve $\chi$\/, 
and hence the so-called adjoint divisors. 
The adjunction divisor of $\chi$ is nothing but the effective divisor 
given by the conductor ideal ${\cal C}_{\chi}$ on $\tilde{\chi}$ 
(notice that $\tilde{\chi}$ is the blowing-up of ${\cal C}_{\chi}\,$). 
It can be computed from the resolution forest as follows. 

Let $q_{1},\ldots,q_{l}$ be the branches of ${\cal T}_{\chi}\,$, and let 
$Q_{1},\ldots,Q_{l}$ be the corresponding points of $\tilde{\chi}$\/, 
by identifying $\tilde{\chi}$ to $\chi_{N}\,$. 

For each vertex $p\in{\cal T}_{\chi}$ set 
$$e_{p}=\sum_{j=1}^{l}e_{p,q_{j}}$$
with the convention that $e_{p,q_{j}}=0$ if the branch $q_{j}$ does not pass 
through the vertex $p$. Then, the adjunction divisor is given by 
$${\cal A}=\sum_{j=1}^{l}\left(\sum_{p\in{\cal T}_{\chi}}
m_{p,q_{j}}(e_{p}-1)\right)\,Q_{j}$$
In the sequel, we will denote 
$d_{Q}\doteq d_{q}\doteq\sum_{p\in{\cal T}_{\chi}}m_{p,q}(e_{p}-1)$\/. 
One has 
$$deg\,{\cal A}=\sum_{p\in{\cal T}_{\chi}}e_{p}(e_{p}-1)deg\,P$$
since $deg\,Q_{j}=deg\,P\cdot[k(Q_{j}):k(P)]$ and 
$e_{p,q_{j}}=m_{p,q_{j}}\cdot[k(Q_{j}):k(P)]$ 
for $p$ in the branch $q_{j}$\/.

Now if we want to give the definition of what an adjoint divisor is, we need 
first some notations. Let $P$ be a closed point of the curve $\chi$ embedded 
in ${\cal S}=\PP^{2}$ and consider the domains $R={\cal O}_{\chi,P}$ 
and ${\cal O}={\cal O}_{{\cal S},P}\;$. Thus, the conductor 
$${\cal C}_{P}={\cal C}_{\overline{R}/R}\doteq\{z\in\overline{R}\;|\;
z\,\overline{R}\subseteq R\}$$
is by definition an ideal in $R$ and $\overline{R}$ at the same time. 
As an ideal of $R$, there exists another ideal $\mbox{\frack A}_{P}$ 
containing the kernel of the natural morphism ${\cal O}\rightarrow R$ 
such that $\mbox{\frack A}_{P}$ is applied onto ${\cal C}_{P}$ by this 
morphism. The ideal $\mbox{\frack A}_{P}$ is called the ideal of 
{\em germs of adjoints} of $\chi$ at $P$ over $\F$. 
In a global situation, the ideal of adjoints $\mbox{\frack A}$ is defined 
as a sheaf of ideals of ${\cal O}_{\cal S}$ over ${\cal S}$ whose stalk at $P$ 
is either $\mbox{\frack A}_{P}$ when $P\in\chi$, or ${\cal O}_{{\cal S},P}$ 
otherwise\footnote{In other words, $\mbox{\frack A}$ is the preimage of the 
conductor sheaf ${\cal C}_{\chi}$ under the natural morphism 
${\cal O}_{\cal S}\rightarrow{\cal O}_{\chi}\;$. }. 
In fact, for $P\in\chi$ one has $\mbox{\frack A}_{P}={\cal O}_{{\cal S},P}$ 
if and only if $P$ is non-singular; hence $\mbox{\frack A}$ 
has a finite support and can be given by the finite set of data 
$\{\mbox{\frack A}_{P}\;\;|\;\;P\in Sing\,(\chi)\}$\/.

On the other hand, with the above notations and following \cite{Casas}, 
for $P\in{\cal S}$ and $h\in{\cal O}_{{\cal S},P}$ with 
$e_{P}(h)\geq e_{p}-1$ given, denote by $H=div(h)$ the divisor defined by 
$h$ on the surface ${\cal S}$\/, and consider 
$\pi_{P}^{\ast}H=div\,(\pi_{P}^{\ast}h)=(e_{p}-1)\,E_{P}+\tilde{H}$\/, 
where $\pi_{P}$ denotes the blowing-up at $P$ 
and $E_{P}$ the exceptional divisor of $\pi_{P}\,$. 
Then $\tilde{H}$ is called the {\em virtual transform} of $H$ 
(with respect to $P$ and the weight $e_{p}$\/), and the multiplicity 
$\mu_{q}(h)\doteq e_{q}(\tilde{H})$ (for $q$ proximate to $p$\/, that is, 
the corresponding point $Q$ is in the strict transform of the 
exceptional divisor created in the blowing-up of the point $P$\/) 
is called the {\em virtual multiplicity} of $h$ at $q$ related to $e_{p}-1$\/. 
By induction, if one substitutes the surface ${\cal S}$ by the 
corresponding one at the inductive step, and by taking the successive 
virtual transforms related to the values $e_{r}-1$, one has in a similar way 
the concept of virtual multiplicity at any $q$ in ${\cal T}_{\chi}$ \/, where 
we take in successive steps the virtual multiplicity $\mu_{r}(h)$ instead of  
the value $e_{p}(h)$ taken in the first step. Then, one has 
$$\mbox{\frack A}_{P}=
\{h\in{\cal O}_{{\cal S},P}\;|\;\mu_{q}(h)\geq e_{q}-1\;\;
\forall q\geq p,\;\;q\in{\cal T}_{\chi}\}$$

As a consequence, for a $\F$\/-rational divisor $D$ on the surface ${\cal S}$ 
one has four equivalent ways to say that $D$ is an adjoint divisor, 
as follows:

\begin{description}

\item[(i)] Adjoint by branches: if the intersection multiplicity of $D$ and 
$\chi$ at every rational branch $q$ of $\chi$ is at least the coefficient 
$d_{q}$ that appears in the adjunction divisor ${\cal A}_{\chi}$\/. 

\item[(ii)] Divisorial adjoint: if ${\bf N}^{\ast}D\geq{\cal A}$\/, where 
${\bf N}=i\circ{\bf n}$\/, ${\bf n}$ being the normalization of $\chi$ 
and $i$ the embedding of $\chi$ in ${\cal S}$\/. 

\item[(iii)] Arithmetic adjoint: if the local equation of $D$ 
at every point $P\in\chi$ is in $\mbox{\frack A}_{P}\,$. 

\item[(iv)] Geometric adjoint: if the virtual multiplicity of $D$ at every 
infinitely near point corresponding to ${\cal T}_{\chi}$ is greater or equal than the 
effective multiplicity of the strict transform of $\chi$ at this point minus one. 

\end{description}

Adjoints are useful to describe the vector space of finite dimension 
$${\cal L}(G)\doteq\{f\in\F(\tilde{\chi})\;|\;(f)+G\geq0\}\cup\{0\}$$
for an arbitrary $\F$\/-rational divisor $G$ on $\tilde{\chi}$\/, 
as derived from the classical Brill-Noether theorem. Assume that 
$\chi$ is given by the homogeneous polynomial 
$F\in\F[{\rm X}_{0},{\rm X}_{1},{\rm X}_{2}]$\/. 
Take a divisor $G$ on $\tilde{\chi}$ that is rational over $\F$ and consider 
a form $H_{0}\in\F[{\rm X}_{0},{\rm X}_{1},{\rm X}_{2}]$ of degree $n$\/, 
with $n\in\N\setminus\{0\}$\/, defined over $\F$\/, 
not divisible by $F$ and satisfying 
$${\bf N}^{\ast}H_{0}\geq G+{\cal A}$$
Then, the {\bf Brill-Noether theorem} states that 
$${\cal L}(G)=\{\frac{h}{h_{0}}\;\;|\;\;
H\in{\cal F}_{n},\;\;H\notin F\cdot\F[{\rm X}_{0},{\rm X}_{1},{\rm X}_{2}]\;\;
{\rm and}\;\;{\bf N}^{\ast}H+G\geq{\bf N}^{\ast}H_{0}\}\cup\{0\}$$
where $h,h_{0}\in\F(\chi)$ denote respectively the rational functions 
$H,H_{0}$ restricted on $\chi$\/, 
and ${\cal F}_{n}\subset\F[{\rm X}_{0},{\rm X}_{1},{\rm X}_{2}]$ denotes 
the set of forms of degree $n$\/.

\vspace{.2cm}

This result allows us to compute a basis of ${\cal L}(G)$ over $\F$ 
by means of the following algorithm, $G$ being an arbitrary rational divisor.

\newpage

\begin{algor}[Brill-Noether algorithm]

$\;$

\begin{center}

For a given $G$\/, define $J=G+{\cal A}$ and $J_{+}=\max\,\{J,0\}$\/. 

\end{center}

\begin{description}

\item[(1)] Take a large enough $n\in\N$ such that there exists 
$H\in{\cal F}_{n}$ not divisible by $F$ with ${\bf N}^{\ast}H\geq J_{+}\,$, 
for instance 
$n>\max\,\left\{m-1,\Frac{m}{2}+\Frac{deg\,J_{+}}{m}-\Frac{3}{2}\right\}$\/, 
$m=deg\,F$ being the degree of $\chi$ (see \cite{HachTh}). 

\item[(2)] Compute a basis over $\F$ of the vector space 
$$V=\{H\in{\cal F}_{n}\;:\;F|H\;\;{\rm or}\;\;{\bf N}^{\ast}H\geq J_{+}\}
\cup\{0\}$$

\item[(3)] Assumed $n\geq m$\/, compute a set of forms of 
${\cal F}_{n}$ giving a basis over $\F$ of the vector space 
$V'=V/W$\/, where $W=\{A\in{\cal F}_{n}\;:\;F|A\}\cup\{0\}$\/. 

\item[(4)] Choose an element $H_{0}\in V\setminus W$ and compute 
the divisor ${\bf N}^{\ast}H_{0}$\/. 

\item[(5)] Compute a set of forms of ${\cal F}_{n}$ being linearly 
independent over $\F$ which generate (modulo $W$\/) the vector space 
of forms $H$ satisfying ${\bf N}^{\ast}H\geq{\cal A}+R$ (or $H=0$), 
where $R\doteq{\bf N}^{\ast}H_{0}-J$\/. 

\item[(6)] If $\{H_{1},\ldots,H_{s}\}$ is the basis obtained in {\bf (5)} 
and for $i=0,1,\ldots,s$ we denote by $h_{i}\in\F(\chi)$ the functions 
$H_{i}$ restricted to $\chi$\/, then 
$${\cal B}=\left\{\frac{h_{1}}{h_{0}},\ldots,\frac{h_{s}}{h_{0}}\right\}$$
is a basis of ${\cal L}(G)$ over $\F$\/. 

\end{description}

\end{algor}

This algorithm also allows us to determine a basis for the space 
$$\Omega(G)\doteq\{\omega\in\Omega(\tilde{\chi})\;|\;(\omega)\geq G\}\cup\{0\}$$
In fact, for any non-zero differential form $\eta$ 
defined over $\F$ denote $K=(\eta)$ the corresponding canonical 
divisor, which is rational over $\F$\/; then one has the 
$\F$\/-isomorphism 
$${\cal L}(K-G)\rightarrow\Omega(G)$$
$$f\mapsto f\eta$$
for any rational divisor $G$\/. Hence, if 
$\{f_{1},\ldots,f_{s}\}$ is a $\F$\/-basis for 
${\cal L}(K-G)$ then the set $\{f_{1}\eta,\ldots,f_{s}\eta\}$ 
is a basis over $\F$ for $\Omega(G)$\/.

Notice that if we want this algorithm to be effective we must solve 
the following related problems: 

\begin{description}

\item[(a)] Compute the adjunction divisor ${\cal A}$ 
for a plane curve $\chi$\/, what can be done from the 
resolution tree at every singular closed point of $\chi$\/. 
Notice that this was already done by means of symbolic 
Hamburger-Noether expressions. 

\item[(b)] Compute the intersection divisor ${\bf N}^{\ast}H$ 
of a homogeneous polynomial $H$ and the curve $\chi$\/, that is, 
the value $\upsilon_{Q}(H)$ at every rational branch $Q$ of $\chi$\/. 
This can be solved by means of the primitive rational parametrizations of 
such branches also given by their corresponding symbolic Hamburger-Noether 
expressions or, more precisely, by {\em lazy evaluation} of these 
parametrizations, (i.e. evaluation step by step whenever necessary for 
the own computation), since the searched values only depend on the first 
terms of such expansions. 

\item[(c)] For a given rational divisor $J$ and a suitable $n\in\N$, 
compute a basis over $\F$ for the vector space 
$$V(J,n)=\{H\in{\cal F}_{n}\;:\;F|H\;\;{\rm or}\;\;{\bf N}^{\ast}H\geq J\}
\cup\{0\}$$
which is the aim of the next section. Note that it also can be done 
by means of the resolution trees and the rational parametrizations of 
$\chi$ computed again from the symbolic Hamburger-Noether expressions. 

\end{description}

\section{Computing bases for ${\cal L}(G)$}

For a given plane curve $\chi$\/, the computation of a basis for 
${\cal L}(G)$\/, $G$ being a rational divisor over $\tilde{\chi}$\/, 
is reduced, by the Brill-Noether theorem, to compute bases for spaces 
of adjoints of a suitable degree $n$\/. We show in this section how 
to impose the required adjunction conditions from the symbolic 
Hamburger-Noether expressions at every rational branch of $\chi$\/, 
by using the classical ideas of Enriques testing passing conditions. 

In practice we know the polynomial  
$F({\rm X}_{0},{\rm X}_{1},{\rm X}_{2})\in
\F[{\rm X}_{0},{\rm X}_{1},{\rm X}_{2}]$ 
defining the absolutely irreducible curve $\chi$ in the projective plane, 
and we have the data of a divisor $G$ that is rational over $\F$\/, 
involving a finite number of rational branches of $\chi$ 
and their corresponding coefficients. 

We first take a value of $n$ such that there exists an adjoint of 
degree $n$ satisfying 
\footnote{We assume $G\geq 0$, but in general one can consider 
the divisor $J_{+}$ instead of $J={\cal A}+G$\/, according to 
the  above notations. }
$${\bf N}^{\ast}H_{0}\geq{\cal A}+G$$
Now computing the residue $R={\bf N}^{\ast}H_{0}-{\cal A}-G$ 
one has to describe the space of homogeneous polynomials $H$ 
of degree $n$ such that ${\bf N}^{\ast}H\geq{\cal A}+R$\/, 
modulo the multiples of $F$\/. 

The problem of finding $H_{0}$ consists just of imposing to 
$H_{0}$ the condition of being an adjoint together with having 
some extra zeros on the divisor $G$\/. On the other hand, in order 
to go on with the Brill-Noether algorithm to describe ${\cal L}(G)$ 
the problem is again the same but taking $R$ instead of $G$\/. 
Thus we have to study the conditions imposed by the inequality 
${\bf N}^{\ast}H\geq{\cal A}+R$ on a homogeneous polynomial $H$ of 
degree $n$\/, $R$ being an arbitrary effective divisor. 

There are two ways to proceed. For the first one, assume that from the 
symbolic Hamburger-Noether expressions we have computed by lazy evaluation 
the primitive rational parametrizations $(X(Z_{r}),Y(Z_{r}))$ given by 
the corresponding Hamburger-Noether expansions at every branch involved 
in the support of the adjunction divisor ${\cal A}$ and $R$\/. 

The {\em Dedekind formula} allows us to find the coefficient $d_{q}$ 
of ${\cal A}$ at the rational branch $q$\/, which is given by 
$$d_{q}=ord_{t}\left(\Frac{f_{Y}(X(t),Y(t))}{X'(t)}\right)
=ord_{t}\left(\Frac{f_{X}(X(t),Y(t))}{Y'(t)}\right)$$
$(X(t),Y(t))$ being a primitive rational parametrization of $q$ 
(notice that either $X'(t)\neq 0$ or $Y'(t)\neq 0$). 
The algorithm to compute the symbolic Hamburger-Noether expressions 
provides us with as many terms of such a parametrization as we need 
to obtain the above orders in $t$\/, by successive substitution 
and lazy evaluation. 

Now we consider the coefficient $r_{q}$ of $R$ at $q$\/, and thus 
the local condition at $q$ imposed to $H$ by the inequality 
${\bf N}^{\ast}H\geq{\cal A}+R$ is given by 
$$ord_{t}h(X(t),Y(t))\geq d_{q}+r_{q}$$
$h$ being the local affine equation of $H$ in terms of the coordinates 
${\rm X},{\rm Y}$ at the corresponding point $P$\/. Again a suitable number 
of steps of the lazy evaluation suffices to describe the first $d_{q}+r_{q}$ 
monomials of the Taylor expansion of $h(X(t),Y(t))$ as a function of 
the indeterminate coefficients of $H$\/, whose vanishing gives the 
required linear conditions, taking all the possible branches $q$ in 
the support of ${\cal A}$ and $R$\/. 

The second way is just the imposition of {\em virtual passing conditions} 
through the infinitely near points of the {\em configuration} of resolution 
$\mbox{\frack C}_{\chi}$ with virtual multiplicities $e_{p}-1$\/, 
what also yields linear conditions on $H$\/. 
The resolution configuration $\mbox{\frack C}_{\chi}$ stands here for 
the set of points $P$ (at the successive blowing-ups) corresponding to 
the vertices $p\in{\cal T}_{\chi}\,$. Notice that from the symbolic 
Hamburger-Noether expressions one can derive not only the total information 
of $\mbox{\frack C}_{\chi}$ but also the information on bigger configurations 
$\mbox{\frack D}$ obtained by adding to $\mbox{\frack C}_{\chi}$ finitely many 
points with multiplicity $1$ at the end of every branch of ${\cal T}_{\chi}\,$. 
Furthermore, the algorithm to compute the symbolic Hamburger-Noether expressions 
gives us also the weights for the resolution tree and local coordinates 
at every infinitely near point, as we have seen in the previous section. 
On the other hand, we say that a divisor $H$ passes (virtually) through a configuration 
$\mbox{\frack D}$ of infinitely near points of $\chi$ with virtual multiplicities 
$\{\mu_{P}\:|\;P\in\mbox{\frack D}\}$ if the virtual multiplicity of $H$ at every 
point $P$ of $\mbox{\frack D}$ (as defined in section {\bf 3}) is greater or equal than $\mu_{P}$\/, 
generalizing the concept of geometric adjoint given in the above section. 

The total number of imposed linear conditions is 
$$\Sum_{P\in\mbox{\frack C}_{\chi}}\Frac{e_{p}\,(e_{p}-1)}{2}deg\,P=
\Frac{1}{2}deg\,{\cal A}$$
since the condition $\mu_{p}(h)\geq e_{p}-1$ is equivalent to the 
vanishing of $\Frac{e_{p}-1}{2}e_{p}$ coefficients, what yields 
this number of conditions over a field 
isomorphic to the residual field $k(P)$\/, and thus 
$\Frac{1}{2}e_{p}\,(e_{p}-1)\,deg\,P$ conditions over the base field $\F$\/. 
Moreover, such conditions are linear independent whenever $n\geq m-3$, 
because of the {\em Noether's adjunction theorem}\/, 
which is refered to the next section, and the virtual transform $\tilde{H}$ 
of $H$ can be computed from the symbolic Hamburger-Noether expressions. 
Note that the first $e_{p}-1$ terms of the Taylor expansion of 
$\tilde{H}(X(t),Y(t))$ vanish. 

Now we must add to ${\bf N}^{\ast}H\geq{\cal A}$ the conditions given by $R$\/. 
If $supp\,R$ does not contain any singular point (that is, the adjoint defined 
by $H_{0}$ passes through $\mbox{\frack C}_{\chi}$ with actual multiplicities 
$e_{p}-1$), then the condition ${\bf N}^{\ast}H\geq{\cal A}+R$ is equivalent 
to ${\bf N}^{\ast}H\geq{\cal A}\;$ and $\;{\bf N}^{\ast}H\geq R$ at the same 
time, and thus the method is just the same as before. This situation can be 
assumed if $n$ is large enough, by a theorem of Serre about the vanishing of 
the cohomology, but in practice the estimate of such values of $n$ is very 
hard and we will give an alternative method to proceed.

Denote by $r_{q}$ the coefficient of $R$ at the rational branch $q$\/, being 
$r_{q}\geq 0$ by assumption. We will show that ${\bf N}^{\ast}H\geq{\cal A}+R$ 
can also be described with virtual passing conditions on $H$\/. In fact, 
consider the configuration $\mbox{\frack C}_{\chi}^{+,R}$ given by 
adding to $\mbox{\frack C}_{\chi}$ the first $r_{q}$ points of 
multiplicity $1$ in the sequence of infinitely near points corresponding 
to the branch $q$\/, for all $q$ in the support of $R$\/.

Recall that the condition ${\bf N}^{\ast}H\geq{\cal A}+R$ can be written 
in terms of the local conditions 
$$ord_{t}h(X_{q}(t),Y_{q}(t))\geq d_{q}+r_{q}\;\;\;\;\;\;(\star)$$
for each rational branch $q$ in $\mbox{\frack C}_{\chi}^{+,R}\;$, 
$(X_{q}(t),Y_{q}(t))$ being a primitive rational parametrization 
corresponding to $q$\/. 
From the inequalities $(\star)$ one gets the following result. 

\vspace{.5cm}

\begin{prop}

Under the above conditions, the inequality ${\bf N}^{\ast}H\geq{\cal A}+R$ 
is equivalent to the condition that the hypersurface defined by $H$ passes 
through the points of $\mbox{\frack C}_{\chi}^{+,R}$ with virtual 
multiplicities $e_{p}-1$ if $p\in\mbox{\frack C}_{\chi}\,$ and $1$ if 
$p\in\mbox{\frack C}_{\chi}^{+,R}\setminus\mbox{\frack C}_{\chi}\,$. 

\end{prop}

\dem

If ${\bf N}^{\ast}H\geq{\cal A}+R$ then ${\bf N}^{\ast}H\geq{\cal A}$\/, 
since $R\geq 0$. Thus, $H$ passes through the points 
$p\in\mbox{\frack C}_{\chi}$ with virtual multiplicities $e_{p}-1$. 
On the other hand, the formula $(\star)$ shows that the virtual transform 
of $H$ at the first point of multiplicity $1$ corresponding to the branch  
$q$ has intersection multiplicity at least 
$r_{q}$ with the strict transform of this branch; 
hence, $H$ passes through the last $r_{q}$ points of 
$\mbox{\frack C}_{\chi}^{+,R}\setminus\mbox{\frack C}_{\chi}$ 
corresponding to $q$ with virtual multiplicity $1$. 

Conversely, if $H$ passes through the points of 
$\mbox{\frack C}_{\chi}^{+,R}$ with the above virtual multiplicities, then 
$(\star)$ is satisfied for any branch $q$ in $\mbox{\frack C}_{\chi}^{+,R}\,$. 

\findemo

\begin{nota}

The above result is considered in \cite{Casas} 
in the case $r_{q}=e_{{\rm N}(q),q}-1$ to study the behaviour of the  
polar curve of a plane curve in characteristic $0$. 
We have proved that in fact the result is also true in any characteristic 
and for arbitrary values of $r_{q}$ whenever $r_{q}\geq 0$. 
Notice that (in total) one considers a number of linear conditions equal to 
$\Frac{1}{2}deg\,{\cal A}+deg\,R$\/, but they may not be linearly independent. 

\end{nota}

\begin{nota}

The theory of Enriques on plane curves with assigned singularities or, 
in more modern terms, the theory of Zariski-Lipman of complete ideals, 
allows us to substitute the weights $e_{p}-1$ in $\mbox{\frack C}_{\chi}\,$ 
and $1$ in $\mbox{\frack C}_{\chi}^{+,R}\setminus\mbox{\frack C}_{\chi}$ 
by other weights $\overline{e_{p}}$ over $\mbox{\frack C}_{\chi}^{+,R}\,$ 
satisfying the so-called proximity inequalities, that is 
$$\overline{e_{p}}\geq\Sum_{r\rightarrow p}\overline{e_{r}}\;\;\;\;\;\;
\forall p\in\mbox{\frack C}_{\chi}^{+,R}$$
This substitution can be done by means of an combinatorial algorithm 
known as the principle of discharge (see for instance \cite{Casas}). 
This algorithm is combinatorial in the sense that one can describe it 
just in terms of the embedded resolution forest associated to the 
configuration $\mbox{\frack C}_{\chi}^{+,R}\,$.

\end{nota}

Also notice that the $r_{q}$ added points of multiplicity $1$ in each branch 
$q$ can be deduced in practice from the symbolic Hamburger-Noether expressions 
computing the first $r_{q}$ terms of the Taylor expansion of the implicit 
function given by the polynomial $g({\rm Z}_{r},{\rm Z}_{r-1})$\/. 
As a consequence of all what we have exposed so far, we can state the 
following result.

\begin{thm}

For any absolutely irreducible plane curve $\chi$ defined over a finite field 
$\F$ and given by a polynomial $F({\rm X}_{0},{\rm X}_{1},{\rm X}_{2})\in
\F[{\rm X}_{0},{\rm X}_{1},{\rm X}_{2}]$\/, and for any rational divisor $G$ 
on $\tilde{\chi}$\/, there exists an algorithm which computes bases over 
$\F$ for ${\cal L}(G)$ consisting of the following steps: 

\begin{description}

\item[(1)] Compute the closed points of the projective plane which are 
singular for the curve $\chi$\/. 

\item[(2)] Compute the symbolic Hamburger-Noether expressions at every 
singular closed point of $\chi$ by using successive symbolic extensions 
of the base field $\F$. 

\item[(3)] Compute an adjoint $H_{0}$ for $\chi$ of degree 
$n\geq m-3$ satisfying ${\bf N}^{\ast}H_{0}\geq{\cal A}+G$\/, 
where ${\cal A}$ is the adjunction divisor of the curve $\chi$ 
computed by means of the step {\bf (2)}, and then compute the 
residue $R={\bf N}^{\ast}H_{0}-{\cal A}-G$\/. 

\item[(4)] Describe the linear conditions ${\bf N}^{\ast}H\geq{\cal A}+R$ 
in terms of the coefficients of a generic form $H$ of degree $n$ and the 
lazy parametrizations of the rational branches of $\chi$ computed also from 
{\bf (2)}, by using either the method given in this section or the 
principle of discharge. 

\item[(5)] Apply the previous steps to get, by using the Brill-Noether 
method, a $\F$\/-basis for the vector space ${\cal L}(G)$\/. 

\end{description}

\end{thm}

\begin{nota}

\begin{description}

\item[] $\;$

\item[i)] The computation of the needed symbolic extensions of $\F$ 
requires factorization of polynomials in one variable, what has an 
effective solution in computational algebra. 

\item[ii)] In fact, we could apply the method to any computable perfect 
field $\F$, that is, when the operations in $\F$ can be done in an effective 
way (for instance, when $\F$ is any field of algebraic numbers). 

\end{description}

\end{nota}

\section{Computing Weierstrass semigroups}

As we will see later, the decoding procedure of Feng and Rao is just based 
on the computation of a basis for ${\cal L}(lP)$\/, $P$ being a rational 
point of $\tilde{\chi}$\/, in the way that if $l\in\Gamma_{P}\,$, the 
Weierstrass semigroup $\Gamma_{P}$ consisting of the Weierstrass non-gaps at 
$P$\/, then such a basis is obtained by adding to a basis of ${\cal L}((l-1)\,P)$ 
a function $f_{l}$ with a unique pole at $P$ of order $l$\/. 
What we are going to do now is to show how one can 
compute the semigroup $\Gamma_{P}$ and the functions $f_{l}$ 
in a quite general situation by using the theory of adjoints. 
For this, we make use of the classical adjunction theorem.

Denote by ${\cal A}_{n}$ the set of adjoints of degree $n$ of the curve $\chi$ 
embedded in $\PP^{2}$ and denote ${\bf N}=i\circ{\bf n}$\/, ${\bf n}$ being 
the normalization of $\chi$ and $i$ the embedding of $\chi$ in $\PP^{2}$\/. 
For every $D\in{\cal A}_{n}$ one can consider 
its {\em pull-back}\/, which is given by ${\bf N}^{\ast}D={\cal A}+D'$ 
for certain $D'\,$. The {\bf adjunction theorem}, due to Noether, says that 
if $n+3\geq deg\,\chi$ the divisors $D'={\bf N}^{\ast}D-{\cal A}$ 
for $D\in{\cal A}_{n}$ are exactly those in the complete linear system 
$|K_{\tilde{\chi}}+(n-m+3)\,L|$\/, $K_{\tilde{\chi}}$ being a canonical 
divisor on $\tilde{\chi}$\/, $L$ the hyperplane section divisor and 
$m=deg\,\chi$ (see \cite{Gor} for details).

This result means that local adjunction conditions are linearly independent 
if imposed on divisors of large enough degree, and this independence is 
in fact global, that is, when imposed on all the points of $\chi$ at 
the same time. In particular, if $n=m-3$ one obtains the following result. 

\begin{prop}

For $n=m-3$ one has an $\F$\/-isomorphism of complete linear systems 
$${\cal A}_{n}\rightarrow|K_{\tilde{\chi}}|$$
$$D\mapsto{\bf N}^{\ast}D-{\cal A}$$

\end{prop}

Notice that this map is injective since $n<m$\/, and the dimension 
over $\F$ of the vector space of forms of dimension $m-3$ in three 
variables equals the arithmetic genus $p_{a}(\chi)$\/. But now the 
total number of linearly independent adjunction conditions is 
$\Frac{1}{2}deg\,{\cal A}$\/, and thus the formula of the geometric genus 
$g(\chi)=p_{a}(\chi)-\Frac{1}{2}deg\,{\cal A}$ can be seen as a problem of 
virtual conditions through the configuration of resolution 
$\mbox{\frack C}_{\chi}\,$.

Under the same hypothesis as in the previous sections, assume that 
$G=lP$\/, where $l$ is a non-negative integer and $P$ is a rational 
point of $\tilde{\chi}$\/, that is, a rational branch defined over $\F$ 
at a certain point of the curve $\chi$\/. 
Then the Riemann-Roch formula can be applied to the divisors $lP\,$ and 
$\,(l-1)\,P$\/, what yields the equality 
$$(\ell(lP)-\ell((l-1)\,P))-(i(lP)-i((l-1)\,P))=1$$
being $0\leq\ell(lP)-\ell((l-1)\,P)\leq 1\,$ 
and $\,-1\leq i(lP)-i((l-1)\,P)\leq 0$. Therefore one has  $l\notin\Gamma_{P}$ 
if and only if $l\geq 1$ and there exists a differential form which is regular 
on $\tilde{\chi}$ and with a zero at $P$ of order $l-1$. 

Notice that $l\in\Gamma_{P}$ if $l\geq 2g$\/. From these remarks one can easily 
prove the following result by using the proposition 5.1 .

\begin{prop}

Let $l\in\Z$ such that $1\leq l\leq 2g-2$. Then: 

\begin{description}

\item[(a)] $l\notin\Gamma_{P}$ if and only if there exists a homogeneous 
polynomial $H_{0}$ of degree $m-3$ with ${\bf N}^{\ast}H_{0}\geq{\cal A}+(l-1)P$ 
such that $P$ is not in the support of the effective divisor 
${\bf N}^{\ast}H_{0}-{\cal A}-(l-1)P$\/. 

\item[(b)] There exists $l'\geq l$ with $l'\notin\Gamma_{P}$ 
if and only if there exists a homogeneous polynomial $H_{0}$ 
of degree $m-3$ such that ${\bf N}^{\ast}H_{0}\geq{\cal A}+(l-1)P$\/. 

\end{description}

\end{prop}

As a consequence, the following result provides us with an effective method 
to do the preprocessing of one-point codes by using plane models for the used 
curve in a quite general situation.

\begin{thm}

Under the same assumptions as above, there exists an algorithm founded 
in the theory of adjoints to 
compute the Weierstrass semigroup $\Gamma_{P}$ together with functions 
$f_{l}$ with a pole at $P$ of order $l$ and regular on 
$\tilde{\chi}\setminus\{P\}$\/, for all $l\in\Gamma_{P}\,$. 

\end{thm}

\dem
\begin{description}

\item[Computing the Weierstrass semigroup: ] $\,$

\vspace{.2cm}

Taking $G=(l-1)P$ instead of the divisor $R$ in proposition 4.1 and using 
the configuration $\mbox{\frack C}_{\chi}^{+,G}$ one can impose the linear 
conditions given by ${\bf N}^{\ast}H\geq{\cal A}+(l-1)P$ on forms $H$ of 
degree $m-3$, which are equivalent to virtual passing conditions through 
$q\in\mbox{\frack C}_{\chi}$ with multiplicities $e_{q}-1$ and through 
$q\in\mbox{\frack C}_{\chi}^{+,G}\setminus\mbox{\frack C}_{\chi}$ with 
multiplicity $1$. 

Then for $l$ increasing from $l=0$ (always in $\Gamma_{P}$\/) 
one imposes successively the linear conditions given by 
${\bf N}^{\ast}H\geq{\cal A}+lP$\/, adding one condition in each step. 
Thus, the added condition given by the new $l$ is linearly independent 
of the previous conditions, by using the proposition 5.1 , if and only if 
$l\notin\Gamma_{P}\,$. All the $g$ gaps of $\Gamma_{P}$, and hence 
the semigroup itself, are computed in at most $2g$ steps. 

\item[Computing the functions $f_{l}\,$: ] $\,$

\vspace{.2cm}

There are two ways to proceed. One way is to compute the functions 
$f_{l}$ for all $l\leq\tilde{l}$\/, $\tilde{l}$ being the largest 
non-gap that is needed in the computations with the considered 
one-point code. The other way is to compute first a generator system 
\footnote{This generator system may be the set of all the primitive 
elements in the Weierstrass semigroup, which is contained in the set 
of the first $g+e$ non-gaps, $e$ being the minimum non-zero element in the semigroup. 
However, it is much better to take an Ap\'{e}ry system since then one can easily compute the Feng-Rao 
distance of the code (see \cite{WSFR}). }
for the Weierstrass semigroup and then give the functions only for 
all $l$ in such a system, $\tilde{l}$ being now the largest generator. 
Anyway, the method described below, which is a suitable application of  
the Brill-Noether algorithm, works in both cases. 

\begin{description}

\item[(i)] Compute a homogeneous polynomial $H_{0}$ 
not divisible by $F$ of large enough degree $n$ 
satisfying ${\bf N}^{\ast}H_{0}\geq{\cal A}+\tilde{l}P$\/, 
and take $l\in\Gamma_{P}$ with $l\leq\tilde{l}$\/. 

\item[(ii)] Denoting ${\bf N}^{\ast}H_{0}={\cal A}+lP+R_{l}$ one has 
$R_{l-1}=R_{l}+P$\/, $R_{l}$ being effective. Thus, for decreasing $l$ 
we can impose the conditions ${\bf N}^{\ast}H\geq{\cal A}+R_{l}$ 
by means of the proposition 4.1 in order to find a homogeneous polynomial 
$H_{l}$ of degree $n$ not divisible by $F$ 
such that ${\bf N}^{\ast}H_{l}\geq{\cal A}+R_{l}$ 
but not satisfying the condition ${\bf N}^{\ast}H_{l}\geq{\cal A}+R_{l-1}\,$. 

\item[(iii)] Thus, the function $f_{l}=H_{l}/H_{0}$ restricted to $\chi$ is 
regular on $\tilde{\chi}\setminus\{P\}$ and has a pole at $P$ of order $l$\/. 

\end{description}

\end{description}

\findemo

\begin{ejplo}

Let $\chi$ be the Klein quartic over $\F_{2}$ given by the equation 
$$F({\rm X},{\rm Y},{\rm Z})=
{\rm X}^{3}{\rm Y}+{\rm Y}^{3}{\rm Z}+{\rm Z}^{3}{\rm X}=0$$
whose adjunction divisor is ${\cal A}=0$, since $\chi$ is non-singular. 
We are going to compute the Weierstrass semigroup at $P=(0:0:1)$ with 
the above method. 

Since $P$ is non-singular one easily obtains by lazy evaluation a local 
pa\-ra\-me\-tri\-za\-tion of $\chi$ at $P$ given by 
$$\left\{\begin{array}{l}
X(t)=t^{3}+t^{10}+\ldots\\Y(t)=t
\end{array}\right.$$
In order to get the gaps of $\Gamma_{P}$ one uses adjoints of degree 
$m-3=1$, whose generic equation 
\footnote{Notice that every plane curve is adjoint to $\chi$\/, 
since ${\cal A}=0$. }
is given by  
$$H({\rm X},{\rm Y},{\rm Z})=a{\rm X}+b{\rm Y}+c{\rm Z}$$
and substituting the first terms of the local parametrization at $P$ 
we get 
$$h(X(t),Y(t))=c+bt+at^{3}+at^{10}+\ldots$$
and proceed as in the above theorem: 

\begin{itemize}

\item $l=1$ is obviously the first gap, since $g=p_{a}(\chi)=3>0$, 
but anyway it can also be checked by the method, since $l=0$ impose 
no condition whereas $l=1$ impose the condition $ord_{t}h(X(t),Y(t))\geq 1$, 
which is equivalent to $c=0$. 

\item For $l=2$, the inequality $ord_{t}h(X(t),Y(t))\geq 2$ is equivalent 
to the conditions $c=b=0$, which are linearly independent of those imposed 
by $l=1$, and thus $l=2$ is a new Weierstrass gap. 

\item If $l=3$, then $ord_{t}h(X(t),Y(t))\geq 3$ is again equivalent to 
$c=b=0$. Therefore the new condition depends on the previous one and one 
has $3\in\Gamma_{P}\,$. 

\item At last, when $l=4$ the condition $ord_{t}h(X(t),Y(t))\geq 4$ is 
equivalent to $c=b=a=0$ and one obtains $l=4$ as the third gap of 
$\Gamma_{P}$ and the procedure ends. 

\end{itemize}

Thus the Weierstrass gaps are $l=1,2,4$ 
and the minimal generator system is then $\{3,5,7\}$ 
\footnote{Notice that this semigroup is not symmetric,  
since the conductor is $C=5<6=2g$\/, and in fact this  
set of generators is the Ap\'{e}ry system related to $\ell=3$. }. 
We are going to compute a function $f_{l}$ for each of these 
three generators also with the method described above.

We apply first the Brill-Noether algorithm to $G=7P$ to obtain a form 
$H_{0}$ of degree $n=4$ not divisible by $F$ such that 
${\bf N}^{\ast}H_{0}\geq J_{+}=J=G=7P$\/. 
That is, taking $H_{0}$ as a generic form of degree $4$ with 
coefficients as variables, the needed condition is equivalent 
to $ord_{t}H_{0}(X(t),Y(t),1)\geq 7$, being $(X(t),Y(t))$ the 
above local parametrization. This can be easily tested with a 
computer and one gets for instance the form 
$H_{0}={\rm X}^{2}{\rm Y}{\rm Z}$\/, which is not divisible by $F$\/. 

Now in order to compute ${\bf N}^{\ast}H_{0}$ we use the symmetry 
of $F$ with respect to the three variables to get local parametrizations 
at the points $Q_{1}=(1:0:0)$ and $Q_{2}=(0:1:0)$\/. Thus, one easily obtains 
$${\bf N}^{\ast}H_{0}=2\,{\bf N}^{\ast}({\rm X})+{\bf N}^{\ast}({\rm Y})+
{\bf N}^{\ast}({\rm Z})=7P+4Q_{1}+5Q_{2}$$

Then, in order to get $f_{7}$ we compute $R_{7}=4Q_{1}+5Q_{2}\,$, and find 
with the above method a form $H_{7}$ of degree $4$ not divisible by $F$ 
such that ${\bf N}^{\ast}H_{7}\geq R_{7}$ but not satisfying 
${\bf N}^{\ast}H_{7}\geq R_{6}=R_{7}+P$\/. This is equivalent to the 
condition ${\bf N}^{\ast}H_{7}\geq R_{7}$ together with the local condition 
at $P$ given by 
$$ord_{t}H_{7}(X(t),Y(t),1)=0$$
obtaining for instance $H_{7}={\rm Z}^{4}$ and hence 
$f_{7}=\Frac{{\rm Z}^{3}}{{\rm X}^{2}{\rm Y}}$\/. 

In a similar way one checks that $H_{5}={\rm Y}^{2}{\rm Z}^{2}$ 
satisfies ${\bf N}^{\ast}H_{5}\geq R_{5}$ but not 
${\bf N}^{\ast}H_{5}\geq R_{4}\,$, obtaining 
$f_{5}=\Frac{{\rm Y}{\rm Z}}{{\rm X}^{2}}$\/, 
and $H_{3}={\rm X}{\rm Y}{\rm Z}^{2}$ satisfies 
${\bf N}^{\ast}H_{3}\geq R_{3}$ but not 
${\bf N}^{\ast}H_{3}\geq R_{2}\,$, obtaining 
$f_{3}=\Frac{\rm Z}{\rm X}$\/. 
In particular, a basis of ${\cal L}(7P)$ over $\F_{2}$ is given by 
$$\{1,\Frac{\rm Z}{\rm X},\Frac{{\rm Y}{\rm Z}}{{\rm X}^{2}},
\Frac{{\rm Z}^{2}}{{\rm X}^{2}},\Frac{{\rm Z}^{3}}{{\rm X}^{2}{\rm Y}}\}$$

\end{ejplo}

\vspace{.2cm}

There is an alternative way to get the functions $f_{l}$ from the 
Brill-Noether algorithm. Assume that a basis $\{h_{1},\ldots,h_{s}\}$ 
of ${\cal L}(\tilde{l}P)$ over $\F$ has been already computed and 
that $\tilde{l}$ is not a gap. We propose a triangulation method 
which works by induction on the dimension $s$ as follows: 

\begin{description}

\item[(1)] By computing first the pole orders $\{-\upsilon_{P}(h_{i})\}$ 
at $P$\/, assume that the functions $\{h_{i}\}$ are ordered in such a 
way that these pole orders are increasing in $i$\/. 

\item[(2)] At least the function $h_{s}$ satisfies 
$-\upsilon_{P}(h_{s})=\tilde{l}$ and we set $f_{\tilde{l}}\doteq h_{s}\,$. 
If any other $h_{j}$ satisfies the same condition, there exists 
a non-zero constant $\lambda_{j}$ in $\F$ such that 
$-\upsilon_{P}(h_{j}-\lambda_{j}h_{s})<\tilde{l}$\/; then we change 
such functions $h_{j}$ by $g_{j}\doteq h_{j}-\lambda_{j}h_{s}$ and 
set $g_{k}\doteq h_{k}$ for all the others. The result now is obviously 
another basis $\{g_{1},\ldots,g_{s}\}$ of ${\cal L}(\tilde{l}P)$ over 
$\F$ but with only one function $g_{s}=f_{\tilde{l}}$ whose pole at $P$ 
has maximum order $\tilde{l}$\/. 

\item[(3)] Since the functions $g_{i}$ are linearly independent over $\F$ 
and $-\upsilon_{P}(g_{i})<\tilde{l}$ for $i<s$\/, one has obtained a 
basis $\{g_{1},\ldots,g_{s-1}\}$ of ${\cal L}(l'P)$ over $\F$, where 
$l'$ denotes the largest non-gap such that $l'<\tilde{l}$\/. But now 
the dimension is $s-1$ and we can continue by induction. 

\end{description}

The above procedure also provides us with a function $f_{l}$ for each 
non-gap $l\leq\tilde{l}$\/. In fact, it can be used to compute the 
Weierstrass semigroup up to an integer $\tilde{l}$\/, since the maximum 
non-gap $l'$ such that $l'\leq\tilde{l}$ is just 
$\max\,\{-\upsilon_{P}(h_{1}),\ldots,-\upsilon_{P}(h_{s})\}$\/, 
in the above notations, and so on by induction.

\section{Effective construction of AG codes}

Let $\tilde{\chi}$ be a non-singular projective algebraic 
curve defined over a finite field $\F$ such that $\tilde{\chi}$ 
is irreducible over $\overline{\F}$\/. 
In order to define the Algebraic Geometry codes, 
take $\F$\/-rational points $P_{1},\ldots,P_{n}$ of the curve and 
a $\F$\/-rational divisor $G$ (which can be assumed effective) 
having disjoint support with $D\doteq P_{1}+\ldots+P_{n}\,$, 
and then consider the well-defined linear maps 
$$\begin{array}{c}
ev_{D}\;:\;{\cal L}(G)\longrightarrow\F^{n}\\
f\mapsto(f(P_{1}),\ldots,f(P_{n}))
\end{array}
\;\;\;\;\;\;{\rm and}\;\;\;\;\;\;
\begin{array}{c}
res_{D}\;:\;\Omega(G-D)\longrightarrow\F^{n}\\
\omega\mapsto(res_{P_{1}}(\omega),\ldots,res_{P_{n}}(\omega))
\end{array}$$

\vspace{.1cm}

\noindent One defines the linear codes 
$$C_{L}\equiv C_{L}(D,G)\doteq Im(ev_{D})\;\;\;\; , \;\;\;\;
C_{\Omega}\equiv C_{\Omega}(D,G)\doteq Im(res_{D})$$

The length of both codes is obviously $n$, and one has 
$(C_{\Omega})=C_{L}^{\perp}$ by the {\em residues theorem}\/. 
On the other hand, given $D$ and $G$ as above there exists 
a differential form $\omega$ such that 
$C_{L}(D,G)=C_{\Omega}(D,D-G+(\omega))$ 
and thus it suffices to deal with the codes of type $C_{\Omega}$\/. 

Denote by $k(C)$ and $d(C)$ the dimension over $\F$ and 
the minimum distance of the linear code $C$ respectively, 
$d(C)$ being the minimum number of non-zero entries 
of a non-zero vector of $C$\/. Goppa estimates for $k(C)$ and 
$d(C)$ are deduced from the Riemann-Roch formula as follows 
(see \cite{TsfVla} for further details). If $2g-2<deg\,G< n$; then 

$$(1)\;\;\left[\begin{array}{ccl}
k(C_{L})&=&deg\,G+1-g\\
d(C_{L})&\geq&n-deg\,G
\end{array}\right.
\;\;\;\;\;
(2)\;\;\left[\begin{array}{ccl}
k(C_{\Omega})&=&n-deg\,G+g-1\\
d(C_{\Omega})&\geq&deg\,G+2-2g
\end{array}\right.$$

The main problem to solve for the construction of such codes consists 
of computing bases for ${\cal L}(G)$ (or $\Omega(G-D)$)\/, finding 
points (rational or not) of the curve and evaluating functions of 
${\cal L}(G)$ at some rational points (or computing residues of 
differential forms in $\Omega(G-D)$ at those points). Thus, with the 
assumption of having a (possibly singular) plane model $\chi$ of the 
curve $\tilde{\chi}$\/, and since the codes of type 
$C_{L}$ and $C_{\Omega}$ are not only dual each other but 
both classes of codes are essentially the same (see \cite{TsfVla}), 
the computational algorithms that are involved in these problems 
will basically be reduced to the following ones: 

\begin{description}

\item[(1)] Find all the closed singular points and all the $\F$\/-rational 
points of $\chi$\/, what can be done by means of Gr\"{o}bner bases computation 
(see \cite{HachTh}). 

\item[(2)] Compute the order of a function at a rational point $P$ 
and evaluate the function at this point when possible, what can be done 
from lazy parametrizations at the rational branch corresponding to $P$\/. 
More precisely, if $\phi=G/H$ is a quotient of homogeneous polynomials of 
the same degree in three variables, and $(X(t),Y(t))$ is the rational 
parametrization obtained from the symbolic Hamburger-Noether expressions 
for the branch given by $P$\/, the order can be computed taking at $P$ the 
corresponding local affine equation $g/h$ of $\phi$ and doing the 
substitution 
$$\Frac{g(X(t),Y(t))}{h(X(t),Y(t))}=
\Frac{a_{r}t^{r}+\ldots}{b_{s}t^{s}+\ldots}$$
obtaining the order $r-s$ by lazy evaluation. Moreover, if $\phi$ is 
well-defined at $P$ (what always happens in the applications to Coding 
Theory), then the evaluation of $\phi$ at $P$ is $0$ if $r>s$\/, 
and $a_{r}/b_{s}$ if $r=s$\/, since we actually work at the point $t=0$. 

\item[(3)] Find a basis for ${\cal L}(G)$ using the Brill-Noether method, 
what has been presented in this paper in an original way by means of 
symbolic Hamburger-Noether expressions and testing virtual passing conditions 
following the Brill-Noether algorithm. 

\end{description}

An interesting case is when $G=mP$, $P$ being an extra rational point 
of $\tilde{\chi}$. In this case the codes $C_{m}\doteq C_{\Omega}(D,mP)$ 
can be decoded by the majority scheme of the Feng and Rao algorithm, 
which is so far the most efficient method for the considered codes 
(see \cite{FR}). 

In order to apply this decoding method, one has to fix 
for every non-negative integer $i$ a function $f_{i}$ in $\F(\tilde{\chi})$ with 
only one pole at $P$ of order $i$ for those values of $i$ for which 
it is possible, i.e. for the integers in the Weierstrass semigroup 
$\Gamma=\Gamma_{P}$ of $\tilde{\chi}$ at $P$. For a received 
word $\mbox{\bf y}=\mbox{\bf c}+\mbox{\bf e}$, where $\mbox{\bf c}\in C_{m}$, 
one can consider the unidimensional and bidimensional syndromes 
given respectively by 

$$s_{i}(\mbox{\bf y})\doteq\sum_{k=1}^{n}e_{k}\,f_{i}(P_{k})
\;\;\;\;\;\;\;\;{\rm and}\;\;\;\;\;\;\;\;
s_{i,j}(\mbox{\bf y})\doteq\sum_{k=1}^{n}e_{k}\,f_{i}(P_{k})\,f_{j}(P_{k})$$
 
\noindent Notice that the set $\{f_{i}\;|\;i\leq m,\;\; i\in\Gamma\}$ 
is actually a basis for ${\cal L}(mP)$ and hence one has 
$$C_{m}=\{\mbox{\bf y}\in\F^{n}\;|\;s_{i}(\mbox{\bf y})=0\;\;
{\rm for}\;\;i\leq m\}$$
thus we can calculate $s_{i}(\mbox{\bf y})$ from the received word 
${\bf y}$ as $s_{i}(\mbox{\bf y})=\sum_{k=1}^{n}y_{k}\,f_{i}(P_{k})$ 
for $i\leq m$\/, and such syndromes are called {\em known}\/. 

In fact, 
it is a known fact that if one has a sufficiently large number of unknown syndromes 
$s_{i,j}({\bf y})$ for $i+j>m$ one could know the emitted word $\mbox{\bf c}$, 
and all above syndromes can be computed by majority voting (see \cite{FR}). 
In practice, the main problem is computing $\Gamma$ and 
the functions $f_{i}$ achieving the values of the semigroup $\Gamma$ 
in order to carry out this decoding algorithm. This is just the other 
problem which has been solved in this paper in a general situation, by 
using the symbolic Hamburger-Noether expressions and the theory of adjoints.

\begin{ejplo}

Let $\chi$ be again the curve given in the example 2.4 by the equation 
$$F({\rm X},{\rm Y},{\rm Z})={\rm X}^{10}+{\rm Y}^{8}{\rm Z}^{2}+
{\rm X}^{3}{\rm Z}^{7}+{\rm Y}{\rm Z}^{9}=0$$
defined over $\F_{2}$ and whose genus is $g=14$. 
This curve has $64$ affine rational points over $\F_{8}$ (namely 
$P_{1},\ldots,P_{64}$\/) and only one point $P=(0:1:0)$ at infinity, 
which is the only singular point of $\chi$ and which was treated in the 
above example. Thus, if one takes an integer $m$ with $26<m<64$, 
one can construct a code $C_{m}=C_{\Omega}(D,mP)$\/, where 
$D=P_{1}+\ldots+P_{64}\,$, whose parameters are $[64,77-m,\geq m-26]$\/. 
For example, if $m=51$ then the dimension is $k=26$ and $C_{51}$ 
corrects any configuration of 12 errors. In order to construct such a 
code, one has first to compute the vector space ${\cal L}(51\,P)$ by 
the Brill-Noether method and then triangulate to obtain the Weierstrass 
semigroup and the functions, since this one-point code can be decoded 
by the Feng-Rao procedure. In order to illustrate the method, 
we consider a smaller space, namely ${\cal L}(15\,P)$\/. 

By using the symbolic Hamburger-Noether expressions computed in the 
example 2.4 one gets, by means of the Brill-Noether algorithm, a 
basis for such a space, namely 
$$\{ h_{1}=1,h_{2}={\rm X}+{\rm X}^{5}+{\rm Y}^{4},
h_{3}={\rm Y}+{\rm Y}^{4}+{\rm X}^{5},h_{4}={\rm X}^{5}+{\rm Y}^{4},
h_{5}={\rm X}{\rm Y}^{4}+{\rm Y}^{16}+{\rm X}^{20} \}$$
The values at $P$ of these functions are $0,12,12,12,13\,$; thus 
$13$ is the largest non-gap in the range $[0,15]$ and $f_{13}=h_{5}$ 
is a function achieving such value. We fix now $f_{12}=h_{4}$ and 
triangulate $g_{2}=h_{2}+h_{4}$ and $g_{3}=h_{3}+h_{4}\,$, according 
to the triangulation method described in the previous section, and one 
finally obtains the values $0,8,10,12,13$ with the corresponding functions. 
In fact, that is enough to construct the whole Weierstrass semigroup and 
all the possibly needed functions, since the sequence $\{8,12,10,13\}$ is 
telescopic and thus generates the semigroup 
(see \cite{KirPel}). 
Finally, by evaluating those functions at the points $P_{1},\ldots,P_{64}$ 
one easily obtains a parity check matrix for the code $C_{m}\;$. 

\end{ejplo}

As a conclusion, our main contribution to the construction of 
AG codes is a new effective solution to the problems which are 
involved in such a construction by using symbolic Hamburger-Noether 
expressions of a plane model for the smooth curve and testing virtual 
passing conditions, on the basis of the Brill-Noether algorithm. 
This way is simpler than the usual method of blowing-ups and Puiseux 
expansions, in the sense that symbolic Hamburger-Noether expressions 
give at the same time the desingularization and the primitive 
rational parametrizations for the branches of the plane curve. 
On the other hand, we have given an effective solution to the 
general problem of computing the Weierstrass semigroup at a rational 
branch $P$ of a singular plane model by using the theory of adjunction, 
together with functions achieving the pole orders in this semigroup, 
what is essential in the construction and decoding problem of 
one-point codes by means of the majority scheme of Feng and Rao.

\end{document}